\documentclass{ifacconf}

\usepackage{graphicx}      
\usepackage{natbib}        
\usepackage{amsmath}
\usepackage{amssymb}
\usepackage{mathtools}
\usepackage{algorithm}
\usepackage{algpseudocode}
\usepackage{cuted}
\usepackage{multicol}
\usepackage{multirow}
\usepackage{color}
\usepackage{caption}
\usepackage{subcaption}
\usepackage{soul}
\usepackage{xcolor,calc}
\usepackage{booktabs}

\begin{document}
\begin{frontmatter}

\title{Parallel KKT Solver in PIQP for Multistage Optimization \thanksref{footnoteinfo}} 

\thanks[footnoteinfo]{This work was supported as a part of NCCR Automation, a National Centre of Competence in Research, funded by the Swiss National Science Foundation (grant number 51NF40\_225155).}

\author[First]{Fenglong Song} 
\author[First]{Roland Schwan} 
\author[First]{Yuwen Chen}
\author[First]{Colin N. Jones}

\address[First]{Automatic Control Laboratory, École Polytechnique Fédérale de Lausanne, Lausanne, 1015 Switzerland (e-mail: \{fenglong.song, roland.schwan, yuwen.chen, colin.jones\}@epfl.ch).}


\begin{abstract}                
This paper presents an efficient parallel Cholesky factorization and triangular solve algorithm for the Karush–Kuhn–Tucker (KKT) systems arising in multistage optimization problems, with a focus on model predictive control and trajectory optimization for racing. The proposed approach directly parallelizes solving the KKT systems with block-tridiagonal–arrow KKT matrices on the linear algebra level arising in interior-point methods. The algorithm is implemented as a new backend of the PIQP solver and released as open source. Numerical experiments on the chain-of-masses benchmarks and a minimum-curvature race line optimization problem demonstrate substantial performance gains compared to other state-of-the-art solvers.

\end{abstract}


\end{frontmatter}

\section{Introduction}
Optimal control problems (OCPs), which arise in applications such as automotive control and autonomous racing, must often be solved in real time, as the sampling periods are on the order of milliseconds. This stringent requirement makes computational efficiency a critical challenge in solver design. Fortunately, OCPs often exhibit sparse and structured formulations, particularly when discretized from continuous-time Ordinary Differential Equations (ODEs) using multiple shooting or direct collocation that results in block tridiagonal Karush-Kuhn-Tucker (KKT) systems. A rich class of solvers that exploits such temporal sparsity (\cite{hpmpc}, \cite{hpipm}, \cite{fatrop}) has been developed and the computational complexity is proportional to the horizon length $N$, i.e. $O(N)$, under mild assumptions. These methods typically involve solving block-sparse KKT systems via sparse Cholesky or condensed Riccati-based factorizations, which are inherently done sequentially.

Recently, several studies have explored parallelism for the temporal sparsity pattern within OCP solvers. A method of $\mathcal{O}(\log N)$ complexity has been proposed in \cite{sarkka_temporal_2023}. The work proposes an associative operator for the unconstrained linear quadratic OCP that can solve the block tridiagonal KKT system with the parallel scan algorithm (\cite{Harris2007ParallelScan}), making logarithmic-time complexity achievable when a sufficient number of parallel computing units are available. However, the method relies on the fact that the optimal control law admits a closed-form solution in the absence of inequality constraints, making it nontrivial to extend the approach to the constrained settings. To address this limitation, \cite{zhang2025parallelbranchmodelpredictive} incorporates inequality constraints via the augmented Lagrangian method (ALM). Nevertheless, as pointed out in \cite{Pougkakiotis2022}, ALM may struggle to achieve high-precision convergence on difficult problems and can be less reliable than interior point method (IPM) based approaches in such cases.


In this work, we introduce parallelization \emph{directly at the linear-algebra level} for solving the KKT systems without altering the outer optimization algorithm. This design accelerates the computation while preserving the numerical robustness and theoretical properties of the underlying optimization method. Specifically, we extend the parallel Cholesky factorization scheme of \cite{cao2002parallel} to efficiently handle block-tridiagonal–arrow KKT matrices that commonly arise in multistage optimization problems. The proposed parallelization strategy is general and can be integrated into a wide range of quadratic programming solvers that exploit the intrinsic OCP structure.

Our main contributions are summarized as follows:
\begin{itemize}
    \item We extend the parallel Cholesky factorization for block-tridiagonal matrices in \cite{cao2002parallel} to support block-tridiagonal-arrow matrices, which allows us to deal with more general problem types.  We also extend the parallelism for the triangular solve to accelerate the forward and backward substitution when solving the KKT systems, enabling thread-safe parallelism without race conditions.
    \item We analyze the computational complexity of the proposed parallel method and derive an optimal strategy for distributing workload across multiple threads to maximize parallel efficiency. In addition, we quantify the theoretically achievable speedups relative to the sequential factorization in \cite{piqp_multistage}.
    \item We provide an open-source implementation of the proposed multi-threaded method to solve the KKT system and integrate it as a new backend of the PIQP solver from \cite{piqp}\footnote{The code will be made available upon publication.}, enabling seamless use for users. We show that PIQP with our multi-threaded KKT system solver outperforms state-of-the-art single-threaded solvers PIQP, HPIPM and Clarabel through numerical experiments.
\end{itemize}

The paper is organized as follows. Section~\ref{sec:formulation} presents the problem formulation and the structure of the KKT system we aim to solve. Section~\ref{sec:solve_kkt_in_parallel} introduces our parallel Cholesky factorization and triangular solve routines, as well as the optimal strategy to distribute computation across multiple threads and the achievable speedups in theory. Finally, implementation details and numerical results are presented in Section~\ref{sec:results}, showcasing the solver’s performance.

\emph{Notation:} We denote the set of real numbers by $\mathbb{R}$, the set of positive integers by $\mathbb{N}_{+}$, the set of $n$-dimensional real-valued vectors by $\mathbb{R}^n$, and the set of $n\times m$-dimensional real-valued matrices by $\mathbb{R}^{n\times m}$. The set of real symmetric matrices of dimension $n$ is denoted by $\mathbb{S}^n$, and the sets of positive semidefinite and positive definite matrices are denoted by $\mathbb{S}^{n}_{+}$ and $\mathbb{S}^{n}_{++}$, respectively. The floor and ceiling operators are denoted by $\lfloor \cdot \rfloor$ and $\lceil \cdot \rceil$, respectively. For symmetric matrices, we use $\star$ to represent the entries (or blocks) in the upper triangular part for brevity. Finally, $D_{*k}$ denotes the collection of all $D_{ik}$ with a fixed index $k$.

\section{Problem Formulation}\label{sec:formulation}
We consider the multistage optimization problem in \cite{piqp_multistage}:
\begin{equation} \label{eq:stage_problem}
\begin{aligned}
    \min_{x, g} \quad & \sum_{i=0}^{N-1} \ell_i(x_i,x_{i+1},g) + \ell_N(x_N,g) \\
    \text {s.t.}\quad & \bar{A}_ix_i+ \bar{B}_ix_{i+1}+ \bar{E}_ig = \bar{b}_i,\hspace{0.5em} i=0,\dots,N-1, \\
    & \bar{C}_ix_i + \bar{D}_ix_{i+1} + \bar{F}_ig \leq \bar{h}_i,\hspace{0.40em} i=0,\dots,N-1, \\
    & \bar{A}_Nx_N + \bar{E}_Ng = \bar{b}_N, \\
    & \bar{D}_Nx_N + \bar{F}_Ng \leq \bar{h}_N,
\end{aligned}
\end{equation}
with coupled stage cost from stage $0$ to $N-1$
\begin{equation*}
    \ell_i(x_i,x_{i+1},g)\! \coloneqq\! \frac{1}{2}\!\begin{bmatrix}
x_i \\
x_{i+1} \\
g
\end{bmatrix}^\top\!\begin{bmatrix}
\bar{Q}_i & \bar{S}_i^\top & \bar{T}_i^\top \\
\bar{S}_i & 0 & 0 \\
\bar{T}_i & 0 & 0
\end{bmatrix}\!\begin{bmatrix}
x_i \\
x_{i+1} \\
g
\end{bmatrix} \!+ \bar{c}_i^\top x_i,
\end{equation*}
and terminal cost at stage $N$
\begin{equation*}
    \ell_N(x_N,g) \coloneqq \frac{1}{2}\begin{bmatrix}
    x_N \\
    g
    \end{bmatrix}^\top\begin{bmatrix}
    \bar{Q}_N & \bar{T}_N^\top \\
    \bar{T}_N & \bar{Q}_g
    \end{bmatrix}\begin{bmatrix}
    x_N \\
    g
    \end{bmatrix}+ \bar{c}_N^\top x_N + \bar{c}_g^\top g,
\end{equation*}
where $x_i \in \mathbb{R}^{n_i}$ are the stage-wise decision variables, 
$g \in \mathbb{R}^{n_g}$ is a global decision variable, and 
$N \in \mathbb{N}$ is the horizon. 
The matrices $\bar{Q}_i \in \mathbb{S}_+^{n_i}$, 
$\bar{S}_i \in \mathbb{R}^{n_{i+1}\times n_i}$, and 
$\bar{T}_i \in \mathbb{R}^{n_g\times n_i}$ together with 
$\bar{c}_i \in \mathbb{R}^{n_i}$ and 
$\bar{c}_g \in \mathbb{R}^{n_g}$ encode the coupled cost. 
The stage-wise variables are also coupled through equality and inequality constraints encoded by 
$\bar{A}_i \in \mathbb{R}^{p_i\times n_i}$, 
$\bar{B}_i \in \mathbb{R}^{p_i\times n_{i+1}}$, 
$\bar{E}_i \in \mathbb{R}^{p_i\times n_g}$, 
$\bar{b}_i \in \mathbb{R}^{p_i}$, 
$\bar{C}_i \in \mathbb{R}^{m_i\times n_i}$, 
$\bar{D}_i \in \mathbb{R}^{m_i\times n_{i+1}}$, and 
$\bar{h}_i \in \mathbb{R}^{m_i\times n_g}$, respectively.

It should be noticed that formulation \eqref{eq:stage_problem} is more general than the classical OCP formulation, as it can describe the inter-stage coupling not only through the system dynamics but also other general inter-stage costs and constraints. Moreover, it can easily capture the structure where global variables are present, e.g., time-optimal MPC and scenario-based robust/stochastic MPC.






Solving the multistage optimization problem (\ref{eq:stage_problem}) via PIQP~\cite{piqp} includes solving the following Karush-Kuhn-Tucker (KKT) system in each iteration:
\begin{equation}\label{eq:kkt_lin_sys}
    \Psi \Delta x = r,
\end{equation}
where 
\begin{equation}\label{eq:Psi}
\Psi \coloneqq \begin{bmatrix}
\Psi_{0,0} &  \Psi_{1,0}^\top &  0 &  \cdots &  \Psi_{g,0}^\top \\[1ex]
\Psi_{1,0} &  \Psi_{1,1} &  \Psi_{1,2}^\top &  \ddots & \Psi_{g,1}^\top \\
0 &  \Psi_{1,2} &  \ddots &   \ddots &   \vdots \\
\vdots &  \ddots &  \ddots &  \Psi_{N,N} &  \Psi_{g,N}^\top \\[1ex]
\Psi_{g,0} &  \Psi_{g,1} &  \cdots &  \Psi_{g,N} &  \Psi_{g,g}
\end{bmatrix},
\end{equation}

with $\Psi_{i,i} \in\mathbb{S}_{++}^{n_i}$, $\Psi_{g,g}\in\mathbb{S}_{++}^{n_g}$, $\Psi_{i+1,i}\in\mathbb{R}^{n_{i+1}\times n_i}$ and $\Psi_{g,i}\in\mathbb{R}^{n_g,n_i}$. The KKT matrix $\Psi$ is symmetric positive definite (SPD) and has a structured block-tri-diagonal-arrow form in (\ref{eq:Psi}). The arrow blocks $\Psi_{g,i}$ and $\Psi_{g,i}^\top$ reflect the coupling between the stage variables $x_i$ and the global variables $g$. While the detailed derivation is omitted here for brevity, interested readers are referred to \cite{piqp_multistage} for the complete algorithmic framework.

Computing the solution of the KKT system (\ref{eq:kkt_lin_sys}) typically constitutes the computational bottleneck in quadratic programming (QP) solvers. 
In the next section, we present a novel parallel algorithm leveraging multi-core CPU architectures to accelerate computation related to solving~\eqref{eq:kkt_lin_sys}.

In the remainder of this paper, we use 1-based indexing, meaning that the first index is denoted by 1 instead of 0.

\section{Solve KKT System in Parallel} \label{sec:solve_kkt_in_parallel}

\begin{figure*}[tbp]
    \begin{subfigure}[t]{0.48\linewidth}
        \raggedright
        \includegraphics[width=\linewidth]{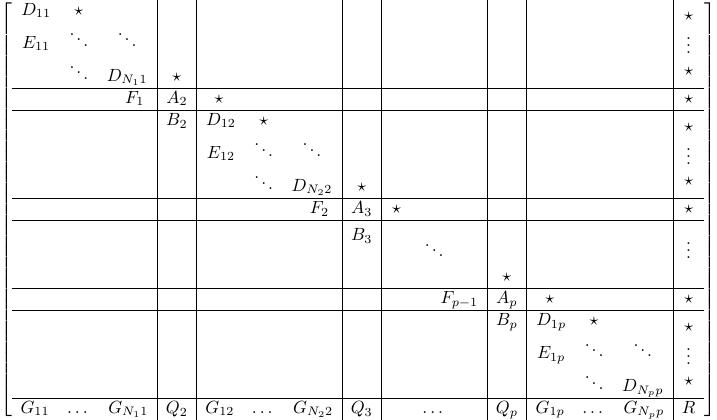}
        \caption{Original KKT matrix $\Psi$.}
        \label{fig:Psi}
    \end{subfigure}
    \hfill
    \begin{subfigure}[t]{0.49\linewidth}
        \raggedleft
        \includegraphics[width=\linewidth]{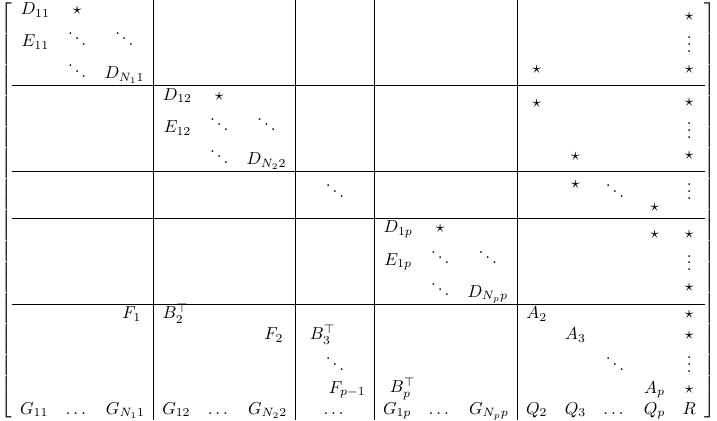}
        \caption{Permuted KKT matrix $\hat{\Psi}\coloneqq P\Psi P$.}
        \label{fig:Psi_hat}
    \end{subfigure}
    \caption{\centering The KKT matrix before and after permutation.}
    \label{fig:Psi_and_Psi_hat}
\end{figure*}

\begin{figure*}[tbp]
    \begin{center}
    \centering
        \includegraphics[width=0.7\linewidth]{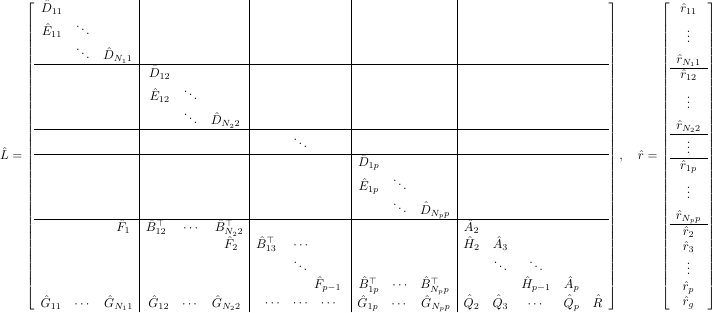}
        \caption{\centering The Cholesky factorization of $\hat{\Psi}$ s.t. $\hat{\Psi} = \hat{L}\hat{L}^\top$ and the permuted right-hand side $\hat{r}\coloneqq P r$.}
        \label{fig:L_hat_and_r_hat}.
    \end{center}
\end{figure*}

\subsection{Parallel Cholesky Factorization}
We summarize our parallel Cholesky factorization in Algorithm \ref{alg:parallel_factor}. The corresponding Basic Linear Algebra Subprograms (BLAS) level operations and flops are marked in the comments in each line. For simplicity, we assume the number of variables at all stages is identical, i.e., $n_i=b$ for $\forall i=0,\ldots, N$. The lines marked with (gv) in comments are only executed when there are global variables $g\in\mathbb{R}^{n_g}$.

The core idea for exposing parallelism in the factorization of the block-tridiagonal-arrow KKT matrix is to break the inter-block dependencies along the time or stage dimension. We view the original KKT matrix $\Psi$ in a different way by organizing the blocks in groups, as illustrated in Figure~\ref{fig:Psi}. The matrix $\Psi$ exhibits strong coupling between neighboring segments of $\{D_{\ast k}, E_{\ast k}, G_{\ast k}\}$ through the coupling blocks $F_k, A_{k+1}, B_{k+1}$, and $Q_{k+1}$. 

To break such coupling, we apply a proper permutation matrix $P$ to obtain $\hat\Psi \coloneqq P \Psi P^\top$ as shown in Figure~\ref{fig:Psi_hat}. With the permutation, we reorder the coupling blocks $\{F_k, A_{k+1}, B_{k+1} \}, \ 1 \le k \le p-1$, to appear in the last rows and columns but before the blocks associated with the global variables. Hence, the processing of the coupling blocks is postponed, allowing the independent portions of the matrix to be factorized first in parallel. This lead to a two-phase routine in Algorithm~\ref{alg:parallel_factor}, including: 
\begin{itemize}
    \item a \emph{parallel} phase, in which operations are distributed among $p$ threads for parallel computation, and
	\item a \emph{sequential} phase, in which the remaining coupled part (bottom-right) is processed serially.
\end{itemize}	

Notice that in Figure~\ref{fig:Psi_hat}, we divide $\hat\Psi$ into \emph{regions} by the solid lines, where each region consists of a group of blocks. This allows us to view the Cholesky factorization of $\hat{\Psi}$ at the higher level in a more compact form. Starting from the top-left, the factorization process proceeds through the following steps. Please notice that in (\ref{eq:region_chol})\textendash(\ref{eq:region_syrk_detailed}), we use~$\hat{~}$ to indicate the final factorization results and~$\tilde{~}$ to indicate intermediate results. 

(1) First, for $k \!=\! 1, \dots, p$, we perform:
\begin{equation}
    {\hat{\Gamma}_k} \coloneqq \begin{bmatrix}
        \hat{D}_{1k} &  &  \\
        \hat{E}_{1k} & \ddots &  \\
         & \ddots & \hat{D}_{N_{k}k} \\
    \end{bmatrix}  = \mathrm{chol}
    \begin{bmatrix}
        {D}_{1k} & \star  &  \\
        {E}_{1k} & \ddots & \ddots  \\
         & \ddots & {D}_{N_{k}k} \\
    \end{bmatrix},
    \label{eq:region_chol}
\end{equation}
\begin{equation}
\hat{\Pi}_k \coloneqq
\begin{bmatrix}
    \hat{B}_{1k}^\top & \cdots & \hat{B}_{N_kk}^\top \\
     &  & \hat{F}_k \\
    \hat{G}_{1k} & \cdots & \hat{G}_{N_kk}
    \end{bmatrix}
    = 
    \begin{bmatrix}
    B_{k}  & &  \\[0.8ex]
     &   & F_k \\[0.8ex]
    {G}_{1k}  & \cdots & {G}_{N_kk} 
    \end{bmatrix} 
    \hat{\Gamma}_k^{-\top},
    \label{eq:region_trsm}
\end{equation}
in which we ignore the operations involving $B_{*k}$ and $\hat{B}_{*k}$ for $k=1$ and those involving $F_k$ and $\hat{F}_k$ for $k=p$. Then update:
\begin{equation}
    \begin{bmatrix}
    \tilde{A}_k & & \star\\
    \tilde{H}_k & \tilde{A}_{k+1} & \star\\
    \tilde{Q}_k & \tilde{Q}_{k+1} & \tilde{R}
    \end{bmatrix} \coloneqq
    \begin{bmatrix}
    {A}_k & &\star\\[0.5ex]
     & {A}_{k+1} &\star \\[0.5ex]
    {Q}_k & {Q}_{k+1}  & {R}
    \end{bmatrix} - \hat{\Pi}_k \hat{\Pi}_k^\top,
    \label{eq:region_syrk}
\end{equation}
where $A_{k}$, $Q_{k}$, $\hat{A}_{k}$, $\hat{Q}_{k}$ are ignored for $k=1$ and $A_{k+1}$, $Q_{k+1}$, $\hat{A}_{k+1}$, $\hat{Q}_{k+1}$ are ignored for $k=p$.

(2) Second, the bottom-right region is factorized as:
\begin{equation}
\hat{\Omega} \coloneqq
     \begin{bmatrix}
    \hat{A}_2 \\
    \hat{H}_2 & \hat{A}_3 & \\
    & \ddots & \ddots & \\
    & & \hat{H}_{p-1} & \hat{A}_p & \\
    \hat{Q}_2 & \hat{Q}_3 & \cdots & \hat{Q}_p & \hat{R}
    \end{bmatrix} =
     \mathrm{chol}\begin{bmatrix}
    \tilde{A}_2 & \star & & & \star \\
    \tilde{H}_2 & \tilde{A}_3 & \ddots & & \star \\
    & \ddots & \ddots & \star & \star\\
    & & \tilde{H}_{p-1} & \tilde{A}_p & \star \\
    \tilde{Q}_2 & \tilde{Q}_3 & \cdots & \tilde{Q}_p & \tilde{R}
    \end{bmatrix}.
\label{eq:region_last_chol}
\end{equation}
The operations in (\ref{eq:region_chol}) and (\ref{eq:region_trsm}) are fully parallelizable since each region $k$ depends only on its own local data. In contrast, (\ref{eq:region_last_chol}) must be executed sequentially for $k = 2, 3, \ldots, p$. The update in (\ref{eq:region_syrk}) contains both parallel and sequential components. Expanding it block-wise gives:
\begin{subequations}\label{eq:region_syrk_detailed}
\begin{align}
& \tilde{A}_k \gets  A_k \!-\! \textstyle\sum_{i=1}^{N_k }\hat{B}_{ik}^{\top} \hat{B}_{ik},  \tilde{Q}_k \gets Q_k \!-\! \textstyle\sum_{i=1}^{N_k} \hat{G}_{ik} \hat{B}_{ik} , \label{eq:tilde_AQ_k} \\
& \tilde{H}_k \gets -\hat{F}_k \hat{B}_{N_k k},  \hspace{1.1cm} \tilde{R}_k \gets -\textstyle\sum_{i=1}^{N_k} \hat{G}_{ik} \hat{G}_{ik}^\top, \label{eq:tilde_HR_k} \\
& \tilde{A}_{k+1} \gets A_{k+1} - \hat{F}_{k} \hat{F}_{k}^{\top},  ~ \tilde{Q}_{k+1} \gets  Q_{k+1} - \hat{G}_{N_kk} \hat{F}_k^\top, \label{eq:tilde_AQ_kp1} \\
& \tilde{R} \gets R + \textstyle\sum_{k=1}^p \tilde{R}_k. \label{eq:tilde_R}
\end{align}
\end{subequations}
Among the operations in (\ref{eq:region_syrk_detailed}), (\ref{eq:tilde_AQ_k}) and (\ref{eq:tilde_HR_k}) can be put into the parallel phase as they involve only data local to region $k$, while (\ref{eq:tilde_AQ_kp1}) and (\ref{eq:tilde_R}) introduce cross-region dependencies and therefore must be put into the sequential phase. This design maximizes parallel efficiency and guarantees thread-safe execution without race conditions.



The sparsity pattern of the lower-triangular factor $\hat{L}$ s.t. $\hat{\Psi} = \hat{L}\hat{L}^\top$ is shown in Figure~\ref{fig:L_hat_and_r_hat}. 
Since factoring $\hat{\Psi}$ introduces fill-ins, i.e., $\{\hat{B}_{2k}^\top, \ldots, \hat{B}_{N_k k}^\top \}$ for $k=2,\ldots,p$ and $\hat{H}_k$ for $k=2,\ldots,p-1$, the total flop count in Algorithm~\ref{alg:parallel_factor} is increased compared to its sequential counterpart~\cite{piqp_multistage}. However, the overall computation time can be reduced thanks to the parallelization.

\begin{algorithm}[t]
\small
    \caption{Parallel Cholesky factorization}
    \begin{algorithmic}[1]
        \Statex \textbf{Parallel Phase}
        \For{$k = 1,...,p$} \textbf{in parallel}
        \State $\hat{R}_k \gets 0$ \Comment{Init}
        \State $\hat{A}_k \gets A_k, \hat{Q}_k \gets Q_k, \hat{B}_{1,k} \gets B_{k}, \hat{B}_{i>1,k} \gets 0$ \textbf{if} $k>1$ \Comment{Init}
            \For{$i=1,...,N_k$}
                \State $\hat{D}_{ik} \gets \text{chol}({D}_{ik}) $  \Comment{potrf, $b^3/3$}
                \State $\hat{E}_{ik} \gets E_{ik} \hat{D}_{ik}^{-\top}$ \textbf{if} $i<N_k$ \Comment{trsm, $b^3$}
                \State $\hat{G}_{ik} \gets {G}_{ik} \hat{D}_{ik}^{-\top} $ \Comment{(gv) trsm, $b_g b^2$}
                \State $\hat{R}_k \gets \! \hat{R}_k \!- \hat{G}_{ik}\hat{G}_{ik}^\top$ \Comment{(gv), syrk, $b_g^2b$}
                \State $\hat{D}_{i+1,k} \gets $ ${D}_{i+1,k}$ $- \hat{E}_{ik} \hat{E}_{ik}^\top$ \textbf{if} $i<N_k$ \Comment{syrk, $b^3$}              
                \State $\hat{G}_{i+1,k} \gets {G}_{i+1,k} - \hat{G}_{ik} \hat{E}_i^\top$ \textbf{if} $i<N_k$ \!\! \Comment{(gv) gemm, $2b_g b^2$}                
                \If{$k>1$}
                    \State $\hat{B}_{i,k}^\top \gets {B}_{i,k}^\top \hat{D}_{ik}^{-\top}$ \Comment{trsm, $b^3$}
                    \State $\hat{A}_k \gets \hat{A}_k - \hat{B}_{ik}^\top \hat{B}_{ik}$  \Comment{syrk, $b^3$}
                    \State $\hat{B}_{i+1,k}^\top \gets {B}_{i+1,k}^\top - \hat{B}_{ik}^\top \hat{E}_{ik}^\top$ \textbf{if} $i<N_k$ \Comment{gemm, $2b^3$}
                    \State $\hat{Q}_k \gets \hat{Q}_k - \hat{G}_{ik} \hat{B}_{ik}$ \Comment{(gv) gemm, $2b_g b^2$}
                \EndIf
            \EndFor
            \State $\hat{F}_k \gets {F}_k \hat{D}_{N_kk}^{-\top}$ \textbf{if} $k<p$ \Comment{trsm, $b^3$}
            
            \State $\hat{H}_k \gets -\hat{F}_k \hat{B}_{N_kk}$ \textbf{if} $1<k<p$  \Comment{gemm, $2b^3$}
            
        \EndFor
        \Statex \textbf{Sequential Phase}
        \For{$k = 2,...,p$}
            \State $\hat{A}_k \gets \hat{A}_k - \hat{F}_{k-1} \hat{F}_{k-1}^\top$ \Comment{syrk, $b^3$}
            \State $\hat{A}_k \gets \text{chol}( \hat{A}_k )$ \Comment{potrf, $b^3/3$}
            \State $\hat{H}_k \gets \hat{H}_k \hat{A}_k^{-\top}$ \Comment{trsm, $b^3$}
            \State $\hat{Q}_k \gets \hat{Q}_k - \hat{G}_{N_{k-1},k-1} \hat{F}_{k-1}^\top$ \Comment{(gv) gemm, $2b_g b^2$}    
            \State $\hat{Q}_k \gets \hat{Q}_k \hat{A}_k^{-\top}$ \Comment{(gv) trsm, $b_g b^2$}
            \State $\hat{A}_{k+1} \gets \hat{A}_{k+1} - \hat{H}_k \hat{H}_k^\top$  \textbf{if} $k<p$  \Comment{syrk, $b^3$}            
            \State $\hat{Q}_{k+1} \gets \hat{Q}_{k+1} - \hat{Q}_k \hat{H}_k^\top$ \textbf{if} $k<p$ \Comment{(gv) gemm, $2b_g b^2$}
            \State $\hat{R} \gets \hat{R} - \hat{Q}_k \hat{Q}_k^\top$  \Comment{(gv) syrk, $b_g^2 b$}
        \EndFor
        \State $\hat{R} \gets \hat{R} + \sum_{k=1}^p \hat{R}_k$ \Comment{(gv) gead $\mathcal{O}(b_g^2)$}
        \State $\hat{R} \gets \text{chol} (\hat{R})$ \Comment{(gv) potrf, $b_g^3/3$}
    \end{algorithmic}
    \label{alg:parallel_factor}
\end{algorithm}

\subsection{Parallel Triangular Solve}

\begin{algorithm}[t]
\small
    \caption{Parallel forward substitution}
    \begin{algorithmic}[1]
    \Statex  \textbf{Parallel Phase}
        \For{$k = 1,..., p$} \textbf{in parallel}
            \State $\Delta \hat{y}_{k} \gets \hat{r}_k$ \textbf{if} $k>1$, $\Delta\hat{y}_{g,k} \gets 0$ \Comment{Initialization}
            \For{$i=1,...,N_k$}
                \State $\Delta \hat{y}_{ik} \gets \hat{D}_{ik}^{-1} \hat{r}_{ik}$  \Comment{trsv, $b^2/2$}
                \State $\Delta \hat{y}_{i+1,k} \gets \hat{r}_{i+1,k} - \hat{E}_{ik} \Delta\hat{y}_{ik}$  ~\textbf{if} $i<N_k$ \Comment{gemv, $2b^2$}
                \State $\Delta \hat{y}_{k} \gets \Delta \hat{y}_{k}  - \hat{B}_{ik}^\top \Delta\hat{y}_{ik}$  ~\textbf{if} $k>1$ \Comment{gemv, $2b^2$}
                \State $\Delta \hat{y}_{g,k} \gets  - \hat{G}_{ik} \Delta \hat{y}_{ik}$ \Comment{(gv) gemv, $2b^2$}
            \EndFor
        \EndFor
    \Statex \textbf{Sequential Phase}
    \For{$k=2,...,p$}
        \State $\Delta\hat{y}_k \gets \Delta\hat{y}_k - \hat{F}_{k-1} \Delta\hat{y}_{N_{k-1},k-1}$  \Comment{gemv, $2b^2$}
        \State $\Delta\hat{y}_k \gets \hat{A}_k^{-1} \Delta\hat{y}_k$  \Comment{trsv, $b^2/2$}
        \State $\Delta\hat{y}_{k+1} \gets \Delta\hat{y}_{k+1} - \hat{H}_k \Delta\hat{y}_{k}$  \textbf{if} $k<p$ \Comment{gemv, $2b^2$}
         \State $\Delta\hat{y}_{g} \gets \Delta\hat{y}_{g} - \hat{Q}_k \Delta\hat{y}_{k}$  \Comment{(gv) gemv, $2b^2$}
    \EndFor
    \State $\Delta\hat{y}_g \gets \hat{r}_g + \sum_{k=1}^p \Delta\hat{y}_{g,k}$ \Comment{(gv) gead, $\mathcal{O}(b_g)$}
    \State $\Delta\hat{y}_g \gets \hat{R}^{-1} \Delta\hat{y}_g$ \Comment{(gv) trsv, $b^2/2$}    
    \end{algorithmic}
    \label{alg:parallel_forward_substitution}
\end{algorithm}

\begin{algorithm}[t]
\small
    \caption{Parallel backward substitution}
    \begin{algorithmic}[1]
    \Statex  \textbf{Sequential Phase}
    \State $\Delta \hat{x}_g \gets \hat{R}^{-\top} \Delta \hat{y}_g$  \Comment{(gv) trsv, $b^2/2$}
    \For{$k=p, ..., 2$}
        \State $\Delta\hat{x}_{k} \gets \Delta\hat{y}_k - \hat{Q}_k^\top \Delta\hat{x}_g$ \Comment{(gv) gemv $2b^2$}
        \State $\Delta\hat{x}_{k} \gets \Delta\hat{x}_{k} - \hat{H}_{k}^\top \Delta\hat{x}_{k+1}$ \textbf{if} $k<p$  \Comment{gemv, $2b^2$}
        \State $\Delta\hat{x}_k \gets \hat{A}_k^{-\top} \Delta\hat{x}_k$  \Comment{trsv, $b^2/2$}
    \EndFor
    \Statex  \textbf{Parallel Phase}
    \For{$k = 1,..., p$} \textbf{in parallel}
        \State $\Delta\hat{x}_{N_kk} \gets \Delta\hat{y}_{N_kk} - \hat{F}_{k}^\top \Delta\hat{x}_{k}$ \textbf{if} $k<p$ \Comment{gemv, $2b^2$}
        \For{$i=N_k,...,1$}
            \State $\Delta\hat{x}_{ik} \gets \Delta\hat{x}_{ik} - \hat{G}_{ik}^\top \Delta\hat{x}_g$  \Comment{(gv) gemv, $2b^2$}
            \State $\Delta\hat{x}_{ik} \gets \Delta\hat{x}_{ik} - \hat{B}_{ik} \Delta\hat{x}_k$ \textbf{if} $k>1$ \Comment{gemv, $2b^2$}
            \State $\Delta\hat{x}_{ik} \gets \hat{D}_{ik}^{-\top} \Delta\hat{x}_{ik}$  \Comment{trsv, $b^2/2$}
            \State $\Delta\hat{x}_{i-1,k} \gets \hat{E}_{i-1,k}^\top \Delta\hat{x}_{ik}$  \textbf{if} $i>1$ \Comment{gemv, $2b^2$}
        \EndFor
    \EndFor
    
    \end{algorithmic}
    \label{alg:parallel_backward_substitution}
\end{algorithm}

After obtaining the factorization of $\hat{\Psi}= \hat{L} \hat{L}^\top$, we now solve the linear system~\eqref{eq:kkt_lin_sys} via the permuted system
\begin{equation*}
    \underbrace{P\Psi P^\top}_{\hat{\Psi}}  \underbrace{P \Delta x}_{\Delta \hat{x}} = \underbrace{P {r}}_{\hat{r}},
\end{equation*}
where $\Delta \hat{x}$ can be computed via a standard forward-backward substitution:
\begin{equation*}
    \hat{L} \Delta\hat{y} = \hat{r}, \quad \hat{L}^\top \Delta \hat{x}  =\Delta\hat{y}.
\end{equation*}
Similar to the factorization, both forward and backward substitutions consist of \emph{parallel} and \emph{sequential} phases, and have been summarized in Algorithm~\ref{alg:parallel_forward_substitution} and Algorithm~\ref{alg:parallel_backward_substitution}, respectively. The subscript notation for $\Delta\hat{x}$ and $\Delta\hat{y}$ is consistent with that of $\hat{r}$, as shown in Figure~\ref{fig:L_hat_and_r_hat}.

\subsubsection{Forward Substitution}
If we view the standard forward substitution on the region-level from the top-left of $\hat{L}$, the process contains the following steps.

(1) First, for each $k = 1, ..., p$, perform:
\begin{subequations}
\begin{align}
    & \begin{bmatrix} \Delta \hat{y}_{1k}^\top & \cdots & \Delta \hat{y}_{N_kk}^\top \end{bmatrix}^\top \gets \hat{\Gamma}_k^{-1} \begin{bmatrix} \hat{r}_{1k}^\top & \cdots & \hat{r}_{N_kk}^\top \end{bmatrix}^\top,  \label{eq:forward_subs_DE}  \\
    & \Delta \tilde{y}_k \gets \hat{r}_k - \begin{bmatrix} \hat{B}_{1k}^\top & \cdots & \hat{B}_{N_kk}^\top \end{bmatrix} \begin{bmatrix} \Delta \hat{y}_{1k}^\top & \cdots & \Delta \hat{y}_{N_kk}^\top \end{bmatrix}^\top, \label{eq:forward_subs_B} \\
    & \Delta \tilde{y}_{g,k} \gets - \begin{bmatrix} \hat{G}_{1k} & \cdots & \hat{G}_{N_kk} \end{bmatrix} \begin{bmatrix} \Delta \hat{y}_{1k}^\top & \cdots & \Delta \hat{y}_{N_kk}^\top \end{bmatrix}^\top, \label{eq:forward_subs_Gk} \\
    & \Delta \tilde{y}_{k+1} \gets \hat{r}_{k+1} - \hat{F}_{N_kk} \Delta x_{N_kk}, \label{eq:forward_subs_F} \\
    & \Delta \tilde{y}_g \gets \hat{r}_g + \textstyle\sum_{k=1}^p \Delta \tilde{y}_{g,k},\label{eq:forward_subs_Gk_sum}
\end{align}
\end{subequations}
where \eqref{eq:forward_subs_B} applied for $k>1$ and \eqref{eq:forward_subs_F} for $k < p$.

(2) Second, for the bottom-right region of $\hat{L}$ and the bottom region of ${r}$:
\begin{equation}
    \begin{bmatrix}
        \Delta \hat{y}_2^\top & \!\! \cdots \! & \Delta \hat{y}_p^\top & \! \Delta \hat{y}_g^\top
    \end{bmatrix}^\top   \!\! \gets 
    \hat{\Omega}^{-1} 
    \begin{bmatrix}
         \Delta \tilde{y}_2^\top & \!\!  \cdots \! & \Delta \tilde{y}_p^\top & \! \Delta \tilde{y}_g^\top
    \end{bmatrix}^\top.
    \label{eq:forward_subs_seq}
\end{equation}

The updates (\ref{eq:forward_subs_DE}), (\ref{eq:forward_subs_B}) and (\ref{eq:forward_subs_Gk}) are put into the parallel phase since they involve only data in the $k$th region, while (\ref{eq:forward_subs_F}), (\ref{eq:forward_subs_Gk_sum}) and (\ref{eq:forward_subs_seq}) are put into the sequential phase. 

\subsubsection{Backward Substitution}
Similarly, we view the steps in backward substitution on the region-level from the bottom-right of $\hat{L}^\top$, which contains the following steps.

(1) First, for the bottom-right region, perform:
\begin{equation}
    \begin{bmatrix}
        \Delta \tilde{x}_2^\top & \!\! \cdots  &\!\! \Delta \tilde{x}_p^\top & \! \Delta \tilde{x}_g^\top
    \end{bmatrix}^\top  \!\! \gets \!
    \hat{\Omega}^{-\top} \!
    \begin{bmatrix}
         \Delta \hat{y}_2^\top &\! \! \cdots &\!\! \Delta \hat{y}_p^\top &\! \Delta \hat{y}_g^\top
    \end{bmatrix}^\top. \!
    \label{eq:backward_subs_seq}
\end{equation}

(2) Second, for each thread $k=p,...,1$:
\begin{subequations}
    \begin{align}
        & \ \Delta \tilde{x}_{N_kk} \gets \Delta \tilde{x}_{N_kk} - \hat{F}_{k+1} \Delta \tilde{x}_{k+1}, \label{eq:backward_subs_F}  \\
        & \begin{bmatrix} \Delta \tilde{x}_{1k} \\ \vdots \\ \Delta \tilde{x}_{N_kk} \end{bmatrix} \gets
        \begin{bmatrix} \Delta y_{1k} \\ \vdots \\ \Delta y_{N_kk} \end{bmatrix} 
        -
        \begin{bmatrix} \hat{G}_{1k}^\top  \\ \vdots \\ \hat{G}_{N_kk}^\top  \end{bmatrix} \Delta\hat{y}_{g} - 
        \begin{bmatrix} \hat{B}_{1k}  \\  \vdots \\ \hat{B}_{N_kk}  \end{bmatrix} \Delta\hat{y}_{k}, \label{eq:backward_subs_BG} \\
        & \begin{bmatrix} \Delta \hat{x}_{1k} \\ \vdots \\ \Delta \tilde{x}_{N_kk} \end{bmatrix}
        \gets 
        \hat{\Gamma}_k^{-1} 
        \begin{bmatrix} \Delta \tilde{x}_{1k} \\ \vdots \\ \Delta \tilde{x}_{N_kk} \end{bmatrix},\label{eq:backward_subs_DE}
    \end{align}
\end{subequations}
where for $k=p$, the $\hat{F}_{k+1}$ term in \eqref{eq:backward_subs_F} must be ignored and similarly the $\hat{B}_{*k}$ term in \eqref{eq:backward_subs_BG} for $k=1$.
The updates (\ref{eq:backward_subs_seq}) and (\ref{eq:backward_subs_F}) must be executed in sequential order, while (\ref{eq:backward_subs_BG}) and (\ref{eq:backward_subs_DE}) can be executed in parallel.

\subsection{Flop Analysis and Optimal Partitioning}

\begin{figure*}[hbtp]
  \begin{subfigure}[t]{0.48\textwidth}
    \raggedright
    \includegraphics[width=\linewidth]{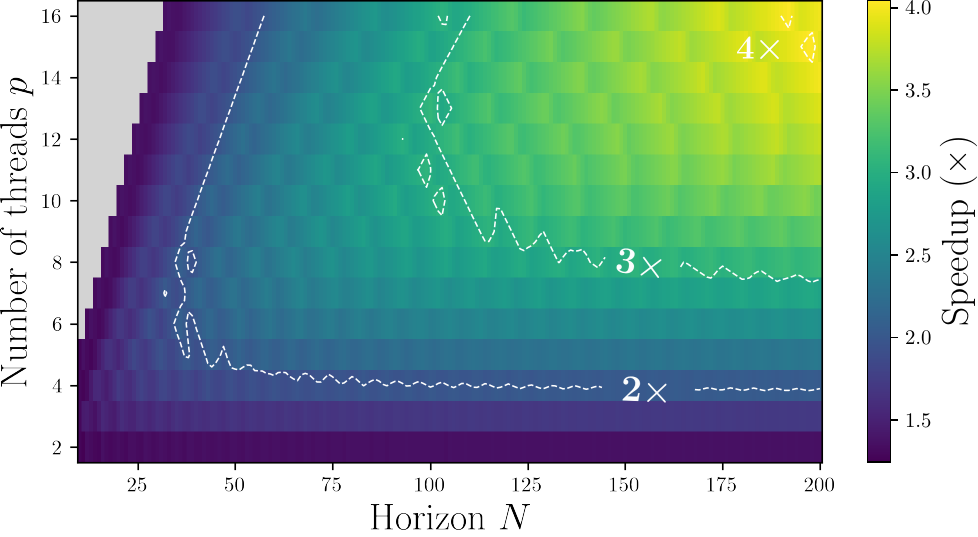}
    \caption{Theoretical Speedup for Cholesky factorization}
    \label{fig:theory_factor_speedup}
  \end{subfigure}
  \hfill
  \begin{subfigure}[t]{0.48\textwidth}
    \raggedleft
    \includegraphics[width=\linewidth]{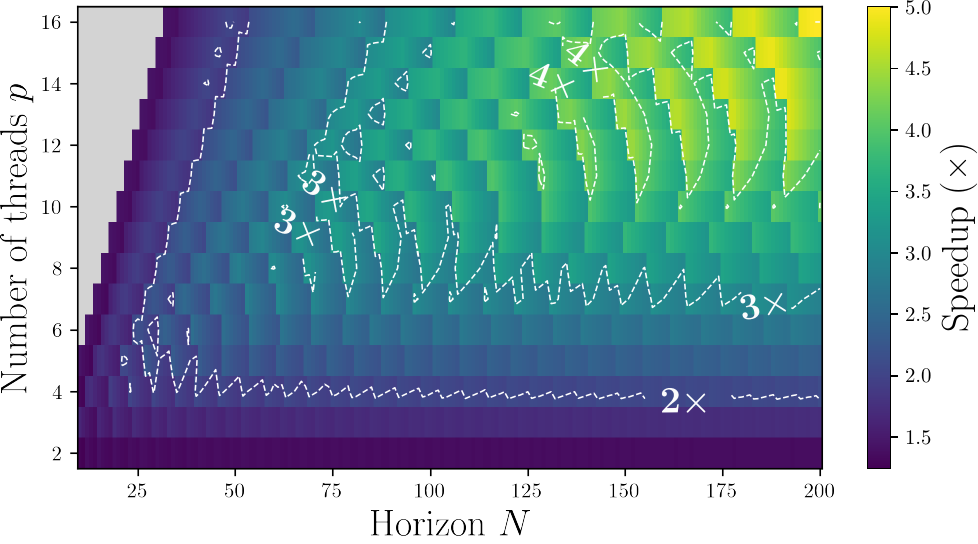}
    \caption{Theoretical Speedup for triangular solve}
    \label{fig:theory_solve_speedup}
  \end{subfigure}
  \caption{Theoretical speedups using our parallel method compared to the sequential method. The gray region indicates the condition $N\geq2p$ is not satisfied, a constraint required for the parallel method across $p$ threads to be meaningful.}
  \label{fig:theory_speedup}
\end{figure*}

In the parallel phase of both the factorization and triangular solve procedures, the first segment does not need to process $B_{*k}$, resulting in a lower computational complexity compared to any $k$-th segment with $k\ge2$. To balance the workload among all threads, the length of the first segment $N_1$ should therefore be larger than that of the remaining segments. For simplicity, we assume that all $n_i$ are equal to $b$ and $n_g=0$ (no global variable), which is a reasonable assumption in nominal OCP formulations.


\begin{table}[h]
\centering
\caption{Flops of KKT solve (w/o global var.) }
\begin{tabular}{@{}lcc@{}}
\toprule
\textbf{Part} & \textbf{Factorization} & \textbf{Triangular Solve} \\
\midrule
Par. phase first seg. & $ (7/3N_1-1)b^3$ & $(5N_1 - 2)b^2$  \\
Par. phase other seg. & $({19}/{3}N_k - 1)b^3$ & $ (9N_k - 2)b^2$ \\
Par. phase last seg. & $({19}/{3}N_p - 4)b^3$ & $ (9N_p - 4)b^2$  \\
\midrule
Seq. phase  & $({10}/{3}p - 16/3)b^3$ & $ (7p-11)b^2$ \\
\bottomrule
\end{tabular}
\label{tab:flops_of_segments}
\end{table}


Given that the Cholesky factorization typically dominates the total computational effort compared to triangular solve, workload balancing across threads requires
\[ (7/3N_1-1)b^3 = (19/3N_k-1)b^3 \Rightarrow \sigma\coloneqq N_1/N_k=19/7. \]
Accordingly, the ideal segment length for the remaining segments is $\bar{N}_k^* \coloneqq  ({N-p+1})/({p+\sigma})$. However, since the segment lengths must be integers, we either round $\bar{N}_k^*$ up to $\lceil \bar{N}_k^* \rceil$ or down to $\lfloor \bar{N}_k^* \rfloor$, depending on which one gives a lower complexity, i.e.,
\begin{equation}
\begin{aligned}
    (N_1^*,N_k^*)= \mathop{\arg\min}\limits_{N_1, N_k\in\mathbb{N}_{+}}  
   & \left\{  \max \left\{ \frac{7}{3}N_1 , \frac{19}{3} N_k \right\}  
    \right\} \\
 \text{s.t.} \hspace{0.3cm} & N_1 + (p-1)N_k + p-1 = N,  \\ 
                            & N_k\in \{ \lceil \bar{N}_k^* \rceil, \lfloor \bar{N}_k^* \rfloor \}.
    \end{aligned}
    \label{eq:partition}
\end{equation}
where $N_1^*$ and $N_k^*$ are the optimal lengths of the first and the other segments.

The time complexity of our method (Algorithm \ref{alg:parallel_factor}) is 
\begin{equation*}
    \mathcal{O}_{par}(N; p) =  \left( \max \left\{ \frac{7}{3} N_1^*, \frac{19}{3}N_k^* \right\} + \frac{10}{3}p- \frac{19}{3} \right) b^3,
\end{equation*}
 compared to the complexity of sequential KKT factorization $\mathcal{O}_{seq} (N)$ given in \cite{piqp_multistage}
\begin{equation*}
    \mathcal{O}_{seq} (N) = \left(\frac{7}{3}N - 2\right) b^3.
\end{equation*}
We illustrate the theoretical speedup of factorization $\gamma := \mathcal{O}_{\mathrm{seq}}(N) / \mathcal{O}_{\mathrm{par}}(N; p)$ in Figure~\ref{fig:theory_speedup} for various numbers of threads $p \in [2, 16]$ and prediction horizons $N \in [10, 200]$. The figure includes both the factorization and triangular-solve stages, following the partitioning strategy described in (\ref{eq:partition}). For a fixed number of threads, the speedup exhibits an overall increasing trend with minor fluctuations as the horizon length $N$ grows. These fluctuations arise from changes in the rounding strategy for $\bar{N}_k^*$, i.e., when the rounding switches between floor and ceiling at certain values of $N$. In addition, the parallel Cholesky factorization and triangular solve must perform additional operations for the fill-in blocks, which offset the gains from parallelization. As a result, the speedup tends to saturate once the horizon $N$ becomes sufficiently long. 

Table~\ref{tab:speedups} outlines the theoretically achievable maximum speedup $\gamma_{\text{max}}$ and the corresponding minimal horizon lengths to achieve $2\times, 3\times, 4\times$ speedups (denoted by $N_{\gamma=2}$ and so forth) and $90\%$ of max speedup (denoted by $N_{0.9\gamma_{\max}}$) with a range of numbers of threads. The results show that with 4 threads, the KKT factorization can achieve a $2.11\times 90\% \approx 1.9$ times speedup once the horizon length $N$ exceeds 43, 
a condition typically met in many practical MPC applications, demonstrating the strong practical applicability of the proposed method.

\begin{table}[htbp]
\centering
\caption{Max speedup in KKT factorization with different number of threads $p$ and min horizon to achieve certain speedups}
\begin{tabular}{@{}c|cccccccc@{}}
\toprule
$p$ & 2 & 4 & 6 & 8 & 10 & 12 & 14 & 16 \\
\midrule
$\gamma_{\text{max}}$ & 1.37 & 2.11 & 2.84 & 3.58 & 4.32 & 5.05 & 5.79 & 6.53 \\
\midrule
$N_{\gamma=2}$  & \textemdash &  83 & 35 & 35 & 41 & 46 & 52 & 58  \\
$N_{\gamma=3}$  & \textemdash & \textemdash & \textemdash & 133 & 101 & 93 & 102 & 102  \\
$N_{\gamma=4}$  & \textemdash & \textemdash & \textemdash & \textemdash & 536 & 244& 201& 190\\
$N_{0.9\gamma_{\max}}$ \!\!\!\!\! &  5 & 43 & 120 & 239 & 384 & 573 & 813 & 1060  \\
\bottomrule
\end{tabular}
\label{tab:speedups}
\end{table}

\section{Numerical Results} \label{sec:results}

We have implemented the proposed factorization and triangular solve routines as a new backend within the solver PIQP (\cite{piqp}). The implementation utilizes C/C++ with the Eigen3 library for the default sparse backend, while leveraging BLASFEO (\cite{blasfeo}) for optimized linear algebra operations in the new backend. Parallelization across multiple threads is achieved using OpenMP. 
In addition to the proposed factorization and triangular solve scheme, we parallelize other naturally parallel operations, including data transfers between Eigen and BLASFEO data structures, as well as block-wise matrix–matrix and matrix–vector multiplications in data preparation. These optimizations further improve throughput and reduce the overall computation time.

We benchmark our algorithm with the sequential multistage solver from PIQP \cite{piqp_multistage} and  HPIPM (\cite{hpipm}). To assess the impact of hardware-level optimizations, our parallel multistage solver as well as sequential PIQP and HPIPM are compiled both with the \texttt{AVX2} instruction set enabled on x86 architecture for fair comparison. All experiments are conducted on an AMD Ryzen 9 3900X processor with 12 cores in a single socket. The turbo boost is disabled to ensure consistent thermal conditions. In addition, we bind OpenMP threads to physical cores to improve performance.

In the remainder of this section, we refer to the PIQP implementation with the general sparse KKT solver in \cite{piqp} as \emph{PIQP (sparse)}, the one with the sequential multistage KKT solver in \cite{piqp_multistage} as \emph{PIQP (seq)}, and to our implementation of the proposed parallel multistage KKT solver as \emph{PIQP (par)}.

\subsection{Chain of Masses Problems}
We evaluate the performance of solvers on the chain-of-masses system from~\cite{chain_mass_problem}, following the same setup as in~\cite{piqp_multistage}. 
A system with $M$ masses includes $2M$ states and $M-1$ inputs. We set $M = 20$ and vary the prediction horizon $N\in \{40, 60 \dots, 200\}$ and the number of threads $p\in\{2,4,\dots,12\}$.

Computation times for the chain-of-masses OCP with varying horizons $N$. The stacked bars show the solver time decomposition into factorization, triangular solve, and other components, distinguished by transparency levels. The arrows indicate the maximum speedup achieved by the fastest PIQP (par) configuration over the PIQP (seq) baseline. For HPIPM, only its total solver runtimes are reported because its internal timings are not exposed. 

\begin{figure*}[htbp]
    \begin{subfigure}[t]{\linewidth}
        \centering
        \includegraphics[width=0.98\linewidth]{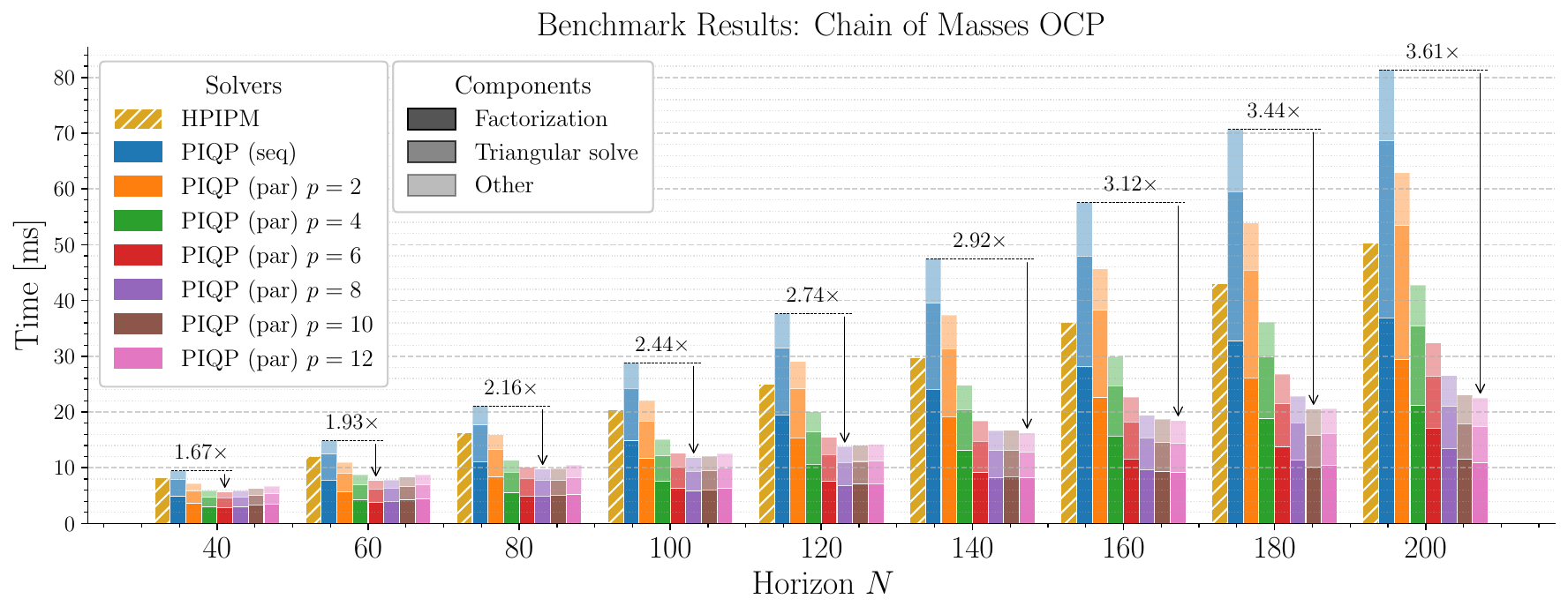}
        \caption{Computation time of PIQP variants and HPIPM for the chain of masses OCP benchmark. The stacked bars show the solver time decomposition into factorization, triangular solve, and other components, distinguished by transparency levels. The arrows indicate the maximum speedup achieved by the fastest PIQP (par) configuration over the PIQP (seq) baseline. For HPIPM, only its total solver runtimes are reported because its internal timings are not exposed. }
        \label{fig:chainmass_time}
    \end{subfigure}
    \begin{subfigure}[t]{0.48\linewidth}
        \centering
        \includegraphics[width=\linewidth]{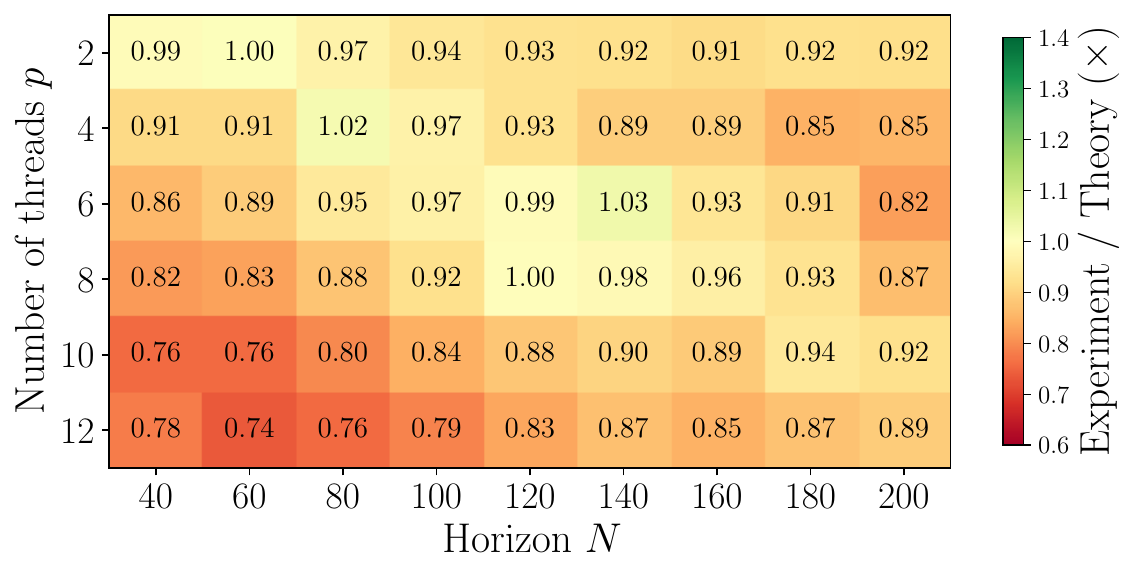}
        \caption{Ratio between experimental and theoretical speedups for the Cholesky factorization stage.}
        \label{fig:speedup_ratio_factor}
    \end{subfigure}
    \hfill
    \begin{subfigure}[t]{0.48\linewidth}
        \centering
        \includegraphics[width=\linewidth]{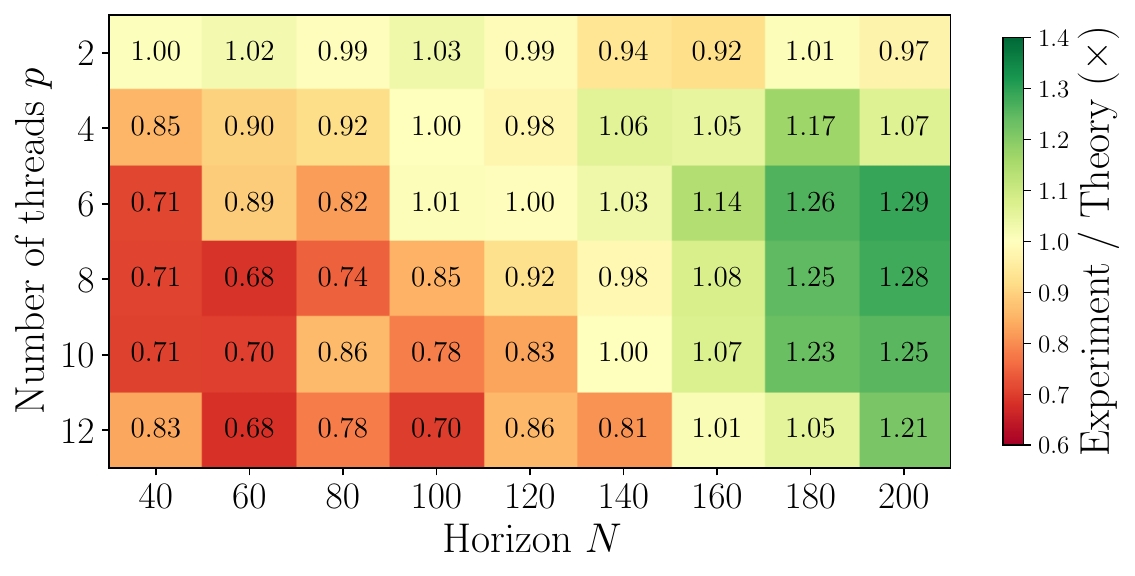}
        \caption{Ratio between experimental and theoretical speedups for the triangular solve stage.}
        \label{fig:speedup_ratio_solve}
    \end{subfigure} 
    \caption{Benchmark results for the chain-of-masses OCP with varying horizons $N$ and numbers of threads $p$. }
    \label{fig:chianmass_benchmark}
\end{figure*}

Figure~\ref{fig:chainmass_time} reports the computation times of multiple variants of PIQP and the state-of-the-art solver HPIPM averaged over 30 runs, as well as the decomposed timings including KKT factorization, triangular solve and other components. 
The results show that the proposed parallel solver achieves substantial acceleration, closely matching the theoretical predictions. For a horizon of $N = 200$ and $p = 12$ threads, PIQP (par) achieves a $3.61\times$ overall speedup over PIQP (seq).
Compared to state-of-the-art solver HPIPM~(\cite{hpipm}), PIQP (par) still achieves up to a $2.24\times$ speedup, demonstrating the effectiveness of the proposed method. 

The experimentally observed speedups generally align well with the theoretical values, as shown in Figure~\ref{fig:speedup_ratio_factor} and Figure~\ref{fig:speedup_ratio_solve}. Slightly lower speedups for $N \le 120$ in both factorization and triangular-solve stages are likely due to the multi-thread scheduling overheads, whereas the experimentally observed triangular-solve speedups for $N \ge 140$ exceed the theoretical predictions, likely due to the parallelization of block-wise matrix–vector multiplications part of the data preparation.

\subsection{Minimum Curvature Race Line Optimization}

We consider the minimum-curvature race line optimization problem, which computes a smooth trajectory for an autonomous race car such that the vehicle can follow it at its handling limits. We consider a quadratic program similar to that in \cite{Heilmeier02102020}, but directly optimize over the spline coefficients of the race line and formulate the problem in a multistage structure.

The race track is represented by a sequence of centerline coordinates $\{x_i^c, y_i^c\}_{i=1}^N$ together with the distances to the boundaries on both left and right sides $\{w_{i}^l, w_{i}^r\}_{i=1}^N$. At each knot, we pre-compute the heading (tangent) direction $t_i \in \mathbb{R}^2$ and right-hand normal direction ${n}_i\in \mathbb{R}^2$. The race line is parametrized by two closed cubic splines in $x$ and $y$, where the $i$th segment is represented by two cubic polynomials respectively:
\begin{equation*}
\begin{aligned}
    x_i(s) & =a_{xi} + b_{xi}s + c_{xi}s^2 + d_{xi}s^3, \\
    y_i(s) &= a_{yi} + b_{yi}s + c_{yi}s^2 + d_{yi}s^3,
\end{aligned}
\end{equation*}
where $s\in[0,1]$ is the normalized curvilinear parameter along the segment. We directly optimize over the spline coefficients $\Theta \coloneqq [\theta_1^\top, \dots, \theta_N^\top]^\top \in \mathbb{R}^{8N}$, where
\[\theta_i \coloneqq[a_{xi}, b_{xi}, c_{xi}, d_{xi}, a_{yi}, b_{i}, c_{yi}, d_{yi}]^\top \in \mathbb{R}^8. \]

\subsubsection{Objective function}
The optimization objective is to minimize the sum of squared curvatures along the whole race line, where for each spline segment it is given by:
\begin{equation*}
\begin{aligned}
    \kappa_i^2 & = \frac{ (x'_i y''_i - x''_i y'_i)^2}{(x_i'^2 + y_i'^2)^{3}}  = \begin{bmatrix}
        x_i'' \\ y_i''
    \end{bmatrix}^\top 
    \begin{bmatrix}
        P_{xx,i} & P_{xy,i} \\ P_{xy,i} & P_{yy,i}
    \end{bmatrix}\begin{bmatrix}
        x_i'' \\ y_i''
    \end{bmatrix} ,
\end{aligned}
\end{equation*}
with 
\begin{equation*}
    P_{xx,i}= \frac{y_i'^2}{z_i}, 
    P_{xy,i}= -\frac{x_i' y_i'}{z_i}, 
    P_{yy,i}= \frac{x_i'^2}{z_i}, 
    z_i = (x_i'^2 + y_i'^2)^3.
\end{equation*}
Since the tangent directions $[x_i', y_i']^\top$ are sufficiently closely aligned with the precomputed centerline heading $t_i$ when the discretization is dense enough, we treat them as constants. With this approximation, the above objective function becomes quadratic w.r.t. $\Theta$.

\subsubsection{Constraints}
We impose continuity constraints up to the second derivative between neighboring spline segments. For simplicity, here we only illustrate for $x$.  Denoting $x_i^\prime(s) \coloneqq \frac{d}{ds}x_i(s)$ and $x_i^{\prime\prime}(s) \coloneqq \frac{d^2}{ds^2}x_i(s)$, the continuity constraints for $i = 1, \dots, N-1$ read:
\begin{equation}
\begin{aligned}
    x_i(1) = x_{i+1}(0),  x_i^\prime(1) = x^\prime_{i+1}(0), x_i^{\prime\prime}(1) = x^{\prime\prime}_{i+1}(0).
\end{aligned}
\label{eq:spline_continuity}
\end{equation}
In addition, the continuity must hold between the first and last segment for a closed race line:
\begin{equation}
\begin{aligned}
    x_N(1) = x_1(0),  x_N^{\prime}(1) = x_1^{\prime}(0),  x_N^{\prime\prime}(1) = x_1^{\prime\prime}(0).
\label{eq:spline_continuity_first_last}
\end{aligned}
\end{equation}
To make sure the race line stays within the track boundary, we impose the following constraints for all knot points: 
\begin{equation}
    {t}_i^\top {r}_i = 0, \quad  -w_i^l \leq {n}_i^\top {r}_i \leq w_i^r
    \label{eq:spline_constraints}
\end{equation}
where ${r}_i := \begin{bmatrix} x_i(0) - x_i^c & y_i(0) - y_i^c \end{bmatrix}^\top $ is the relative position of the $i$th knot point w.r.t. the $i$th centerline point. 

\subsubsection{QP Formulation}
Collecting the objective function and constraints, the minimum curvature race line problem is formulated as a multistage QP with $n_i = 8$ for all stages:
\begin{equation}\label{eq:qp_raceline}
\begin{aligned}
    \min_{\Theta}  ~ \sum_{i=1}^N \kappa_i^2  \quad
    \text{s.t.} ~ (\text{\ref{eq:spline_continuity}}), (\ref{eq:spline_continuity_first_last}), (\ref{eq:spline_constraints}).
\end{aligned}
\end{equation}
Note that the closure constraints (\ref{eq:spline_continuity_first_last}) are equivalently represented by introducing global variables $g$ with $n_g=8$ and enforcing equality constraints between $g$ and $\theta_1$ as well as between $g$ and $\theta_N$. The presence of the couplings leads to an arrow-shaped KKT matrix, as opposed to the block-tridiagonal structure in standard OCP, and thus cannot be efficiently handled by solvers such as HPIPM.

\subsubsection{Benchmark Results}
\begin{table*}[htbp]
\centering
\caption{Runtimes (mean $\pm$ std, in ms) and speedups relative to PIQP (seq) for total solver run, factorization, triangular solve, and other components in the race line optimization problem. Bold values indicate the fastest result within each category.}
\begin{tabular}{@{}l|cc|cc|cc|cc@{}}
\toprule
\multirow{2}{*}{Solver} 
& \multicolumn{2}{c|}{Total} 
& \multicolumn{2}{c|}{Factorization} 
& \multicolumn{2}{c|}{Triangular solve} 
& \multicolumn{2}{c}{Other} \\
& Time [ms] & Speedup & Time [ms] & Speedup & Time [ms] & Speedup & Time [ms] & Speedup \\
\midrule
PIQP (seq) & $137.92 \pm 1.43$ & $1.0\times$ & $51.81 \pm 0.31$ & $1.0\times$ & $63.59 \pm 1.07$ & $1.0\times$ & $22.52 \pm 0.24$ & $1.0\times$ \\  %
PIQP (par) $p=2$ & $98.99 \pm 0.54$  & $1.39\times$ & $38.80 \pm 0.17$ & $1.34\times$ & $43.77 \pm 0.31$ & $1.45\times$ & $16.43 \pm 0.14$ & $1.37\times$ \\ %
PIQP (par) $p=4$ & $68.62 \pm 0.86$  & $2.01\times$ & $26.29 \pm 0.25$ & $1.97\times$ & $28.72 \pm 0.50$ & $2.22\times$ & $13.62 \pm 0.21$ & $1.65\times$ \\ %
PIQP (par) $p=6$ & $56.48 \pm 0.48$  & $2.44\times$ & \textbf{20.40 $\pm$ 0.30} & \textbf{2.54$\times$} & $23.00 \pm 0.26$ & $2.77\times$ & $13.07 \pm 0.13$ & $1.72\times$ \\ %
PIQP (par) $p=8$ & $53.00 \pm 0.57$  & $2.60\times$ & $20.43 \pm 0.34$ & $2.54\times$ & $20.30 \pm 0.29$ & $3.13\times$ & $12.27 \pm 0.13$ & $1.84\times$ \\ %
PIQP (par) $p=10$ & $50.92 \pm 0.47$ & $2.71\times$ & $20.79 \pm 0.48$ & $2.49\times$ & $18.26 \pm 0.04$ & $3.48\times$ & $11.87 \pm 0.05$ & $1.90\times$ \\
PIQP (par) $p=12$ & \textbf{50.30 $\pm$ 0.53} & \textbf{2.74$\times$} & $21.10 \pm 0.37$ & $2.46\times$ & \textbf{17.43 $\pm$ 0.17} & \textbf{3.65$\times$} & $11.77 \pm 0.23$ & $1.91\times$ \\
PIQP (sparse) & $111.40 \pm 1.32$ & $1.24\times$ & $59.32 \pm 0.32$ & $0.87\times$ & $44.14 \pm 0.76$ & $1.44\times$ & \textbf{7.94 $\pm$ 0.22} & \textbf{2.84$\times$} \\
Clarabel & $252.77 \pm 1.58$ & $0.55\times$ & \textemdash & \textemdash & \textemdash & \textemdash & \textemdash & \textemdash \\
\bottomrule
\end{tabular}
\label{tab:speedups_raceline}
\end{table*}

We apply our method to compute the minimum curvature race line for the Silverstone Formula One race track in England that is approximately 5.89km long and divided into 2356 segments, leading to a horizon of $N=2356$ in (\ref{eq:qp_raceline}). Table~\ref{tab:speedups_raceline} summarizes the computation times of PIQP (par) under different numbers of threads, the single-threaded solver PIQP (seq), as well as general-purpose sparse problem solvers PIQP (sparse) and Clarabel (\cite{Clarabel_2024}). The table reports the mean and standard deviation of the total solver time, as well as the contributions from the KKT factorization, triangular solve, and other components if available. As $p$ increases, both the factorization and triangular-solve stages exhibit substantial speedups, leading to an overall improvement of up to $2.74\times$ at $p=12$. The factorization and triangular solve achieve $2.46\times$ and $3.65\times$ speedup, respectively, compared to the sequential baseline.

Interestingly, PIQP (sparse) slightly outperforms the multistage solver PIQP (seq) despite the latter’s BLASFEO-optimized backend. This can be attributed to the extremely sparse structure of the raceline problem, i.e., small blocks compared to the horizon ($N=2356$ and $n_i=8$ per stage) and the blocks themselves being quite sparse. Although the sparse LDL factorization involves less cache-friendly memory access and is marginally slower than the block-dense BLASFEO-based factorization, its triangular solves and residual assembly are significantly cheaper. In contrast, PIQP (seq) performs numerous small dense $8\times8$ and $8\times1$ operations, leading to higher memory overhead.

\section{Conclusion and future work}
We have presented a parallel Cholesky method for solving KKT systems with block–tridiagonal–arrow structures in multistage optimization. By combining a permutation-based decoupling of temporal dependencies with parallel Cholesky factorization and triangular solve, the proposed approach achieves significant runtime reductions while maintaining numerical stability. The method has been integrated into PIQP as a multi-threaded backend and demonstrated on representative benchmarks.

Future work will be extending the approach to a GPU implementation, which allows for $\mathcal{O}(\log N)$ scalability under sufficient parallel computing resources.

\section*{DECLARATION OF GENERATIVE AI AND AI-ASSISTED TECHNOLOGIES}
During the preparation of this work the authors used ChatGPT in order to refine the language, grammar, and readability. After using this tool/service, the author(s) reviewed and edited the content as needed and take(s) full responsibility for the content of the publication.

\bibliography{ifacconf}             
                                                   
\end{document}